\documentclass[12pt,a4paper]{article}

\usepackage{RGpaperSettings} 
\usepackage{multicol}

\title{RG Smoothing Algorithm Which Makes Data Compression}
\author{Anna Sinelnikova}
\date{May 2018}
\affilation{Uppsala University, Sweden}
\begin{document}
\maketitle

\setlength{\columnsep}{20pt} 
\begin{multicols}{2} 
\lettrine[lines=3,slope=-0pt, loversize=0.15]{\color{seccolor}{\ttfamily I}}{would like}
to present a new method of smoothing of one-dimensional curves in arbitrary number of dimensions.
The basic idea is borrowed from \hl{renormalization group (RG) theory} which was applied to biological macromolecules \cite{my}.
In general RG theory is used for \hl{rescaling} of a system. \hlb{You combine the elements into blocks and then treat the blocks the same way as you did with the elements.}
There are two important features that follow from this, that are inherited by our smoothing algorithm:
\begin{itemize}[label=\textcolor{seccolor}{\textbullet}]
\item reduction of the number of elements
\item recursive implementation
\end{itemize}

\noindent These two points distinguish RG smoothing from other methods and make it unique.
 
The first feature means that the smoothing reduces the amount of data.
This property seems very natural if we think about what actually any smoothing method does.
Obviously, it makes the curve or surface smoother.
In a more formal language we can say that smoothing decreases differential curvature of the curve or surface.

But let's now think about the problem from a higher level: form a perspective of the amount of information.
There is some information which we call \hl{data} from where we want to extract the information about the ``ideal data" which can be for example a signal cleaned up from noise.

In any case the point is that the ``useless" information is removed and only the important one remains. So \hlb{any smoothing leads to the loss of at least some information}. The question is if it actually leads to the reduction of the amount of binary data. In our algorithm it does.
 

The feature of the recursive implementation makes the algorithm simpler.
There is only a small instruction which should be repeated over and over again.
In the simplest form, which I will describe in this article, the RG method has the only one tuning parameter: number of recursive steps.
This gives great opportunities for building hardware for the algorithm, because each recursive loop is not complicated and strictly defined.

\section*{Method}\label{sec:method}
RG smoothing algorithm consist of two crucial steps:

\begin{enumerate}[label=\textcolor{seccolor}\theenumi]
\item construction of chain from the data points
\item recursive rescaling procedure for that chain
\end{enumerate}

The RG theory is used in the second step.
So I will start from explaining what is the rescaling procedure for a given chain.
\subsection*{Rescaling procedure}
\begin{figure*}
	\hspace{-1.cm}
	\includegraphics[trim = 0mm 15cm 0mm 0mm, clip,width=1.1\textwidth]{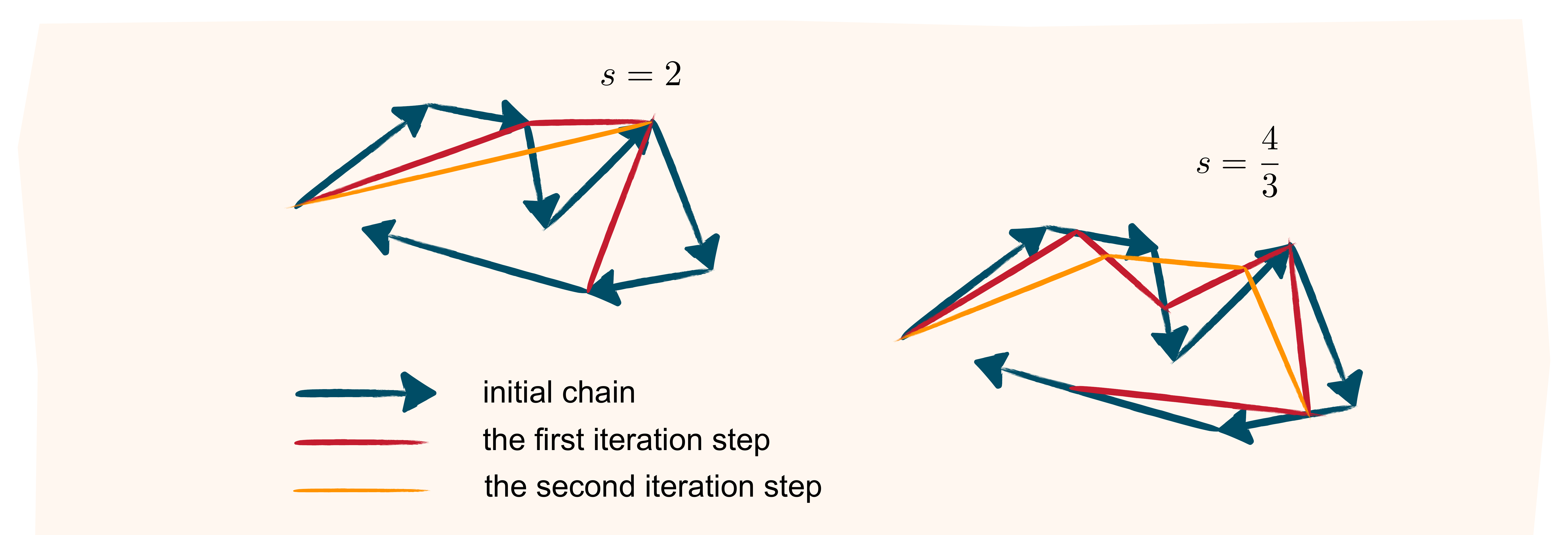}
	\caption{Demonstration of the renormalization procedure. The dark bold arrows connects the original data points. The red thinner lines are result of the first rescaling step. The thinnest yellow line is the result of the second iteration step. Left figure: the scaling parameter is equal to 2; Right figure: the scaling parameter is equal to 4/3.}
	\label{fig:scaling}
\end{figure*}
We assume that we have a chain of rigid segments and want to apply RG theory to it.
As that theory says to us, \hl{we should combine elements into blocks}.
Exactly this process I will call \hl{rescaling}.

The simplest way of rescaling of the chain one can think about is to connect every other site as is
shown on the left part of Figure~\ref{fig:scaling} with red line, where it is applied to the original chain which consists of dark bold arrows.
We could also go for every third site or even more if our chain would be long enough.
The parameter which is responsible for this choice I will call the \hl{scaling parameter $s$}.
More strict it can be defined like this :  \hlb{scaling parameter $s$ is equal to number of old segments combined into a new one}. In the case which is shown in Figure~\ref{fig:scaling}, $s=2$: each new segment connects two old segments.

The second part of the ``definition" for RG says that \hl{we should treat the block in the same way as the elements}.
And this is what \hl{recursive} stands for.
Now we should treat the new chain in the same way as the original one and rescale it again.
This is how we come to the yellow chain (or in our case just one segment) in Figure~\ref{fig:scaling}.

You can see that the number of elements is reduced during the procedure: from 7 segments in the original chain we come to one segment of the yellow chain, through 3 segments of the red chain. We could continue the rescaling procedure if the original chain were longer.



Here it is important to notice that we will loose one segment of the chain every time when number of segments is odd.
During the whole procedure this will unavoidably happen if the initial number of monomers was not a power of 2, as one can easily check.

The RG theory is a complicated mathematical concept with a lot of nuances, which we cannot mention in this article. 
What is important for us is that the RG theory for chains requires the scaling parameter $s$ to be less than~2, or even very close to 1.
The details one can find in the book \cite{grosberg}.

In other words the scaling parameter should lie in the interval $1<s \leq 2$, and thus cannot be an integer number any longer.
The case of non-integer scaling parameter is complicated and allows more than one solution.

I will use one of those solutions which I think is the best in our case.
It is presented in the right plot of Figure \ref{fig:scaling} for scaling parameter $s = 1+1/3$.
The idea is that every time we should connect $s$ old segments.
When $s$ was equal to 2 we connected 2 old segments which means that the new chain consists of each second vertex of old chain.
Now with non-integer $s$ we will not always hit a vertex, but sometimes hit segment itself.
To explain this idea better I will use vectors.

Let's denote every original bold dark arrow in Figure~\ref{fig:scaling}  as $\vec{t}_1$, $\vec{t}_1$, ..., $\vec{t}_7$.
At the first step of rescaling we want to obtain red chain which we will mark as $\vec{t}_1^1$, $\vec{t}_1^1$, ...
So, for our chose of $s=1+1/3$, we can write down according to it's definition :
\setlength{\arraycolsep}{0.0em}
\begin{eqnarray}
\vec{t}_1^1&{}={}&\vec{t}_1 + \frac{1}{3}\vec{t}_2\nonumber\\
\vec{t}_2^1&{}={}&\frac{2}{3}\vec{t}_2 + \frac{2}{3}\vec{t}_3\\
\vec{t}_3^1&{}={}&\frac{1}{3}\vec{t}_3 + \vec{t}_4\nonumber\\
...\nonumber
\label{eq:rescaling1}
\end{eqnarray}
\setlength{\arraycolsep}{5pt}
And  the same for the second recursive step, but with respect to the new red chain ):
\setlength{\arraycolsep}{0.0em}
\begin{eqnarray}
\vec{t}_1^2&{}={}&\vec{t}_1^1 + \frac{1}{3}\vec{t}_2^1\nonumber\\
\vec{t}_2^2&{}={}&\frac{2}{3}\vec{t}_2^1 + \frac{2}{3}\vec{t}_3^1\\
\vec{t}_3^2&{}={}&\frac{1}{3}\vec{t}_3 ^1+ \vec{t}_4^1\nonumber\\
...\nonumber
\label{eq:rescaling2}
\end{eqnarray}
The new chain is marked with yellow color in the same Figure~\ref{fig:scaling}.
The maximum number of scaling steps in this case is 3, as one can see it the plot.

Notice, that we loose again the right end, even after the first step.

There are five important points we should draw into conclusions:
\begin{enumerate}[label=\textcolor{seccolor}\theenumi]
	\item the procedure works in Euclidean space with arbitrary number of dimensioned .
	\item the process is recursive, i.e. it consists of repeated steps.;
	\item in each step we reduce the number of data points;
	\item we can loose the endings of a chain during the procedure;
	\item the procedure does not treat the ends of the chain in the same way.
\end{enumerate}

The first three items in the list are considered as advantages.
We have never mentioned the first advantage before, but indeed, we have never refered to the number of dimensions of the space we consider.
We worked in plane, i.e with 2D space, but obviously the same algorithm will work for 3D space. And since the equations~1 and~2 operate only with vectors, then the only requirement for the space is to have positively defined norm.

Number (4) and number (5) in the list are downsides.
Loosing ends obviously leads to loss of information, which can be crucial for some applications.
Number (5) is a bad thing, because the direction of moving along the chain was introduced artificially, there should be no distinction between two directions.

We will fix both disadvantages in the final algorithm which I am going to present in the next section.
\subsection*{Constructing the Chain}
Now let's come back to the original problem of smoothing.
Assume that we have a set of data points which we want to operate with. It can be signal with some noise or a trajectory in two- or three-dimensional space.
To use the method I have just described, first we should construct a chain from our data.
This is equivalent to enumerate all the data points and connect them in increasing or decreasing order. 
The points can be enumerated in different ways, it is important that the algorithm we present, will work even if there are self crossings and loops in the final chain. 
However, it can affect the efficiency of smoothing.

\section*{Algorithm}\label{sec:algorithm}
Basically we do exactly what is written in the beginning of the previous section, but with an additional trick.
Let me do a small remark about the notation we use: now I will call the recursive procedure as an \hl{iterative} one, because during the smoothing we are approaching the solution more and more accurate with every step.

As we know now, during the procedure we can loose the last data points and thus the information about the end of the data set.
We want to avoid this.

As one could have noted from Figure \ref{fig:scaling} for $s=4/3$  we will not loose any data in the first step if the length of the chain would be $N=4$ instead of $N=7$.
It happens because:
\begin{equation*}
\frac{N}{s} = \frac{4}{\frac{4}{3}} = 3 =\; \rm{integer\;number}.
\label{eq:ints}
\end{equation*}
Then on the one hand we want $s$ always to satisfy this condition which will lead to dependence of $s$ on iteration step.
And on the other hand we have RG requirement for $s$ to be close to 1.

Thus we will require new chain to be one segment shorter that the old one. In other words, let's decrease the number of data points by one in each iteration step. 
This gives the equation for scaling parameter $s$ as a function of iteration step number $p$:
\begin{equation}
s_p^{opt} = \frac{N_p}{N_{p+1}} = \frac{N_p}{N_p -1},
\label{eq:s}
\end{equation}
where $N_p$ is the number of segments at the $p$-th iteration step.
Now we found the smallest $s$ which provides the integer number of segments in a new chain without rounding.

By finding the scaling parameter that way we solved both problems from the previous section as:
\begin{enumerate}[label=\textcolor{seccolor}\textbullet]
	\item no cut off of data set;
	\item the procedure treats the beginning and the end of the data set in the same way.
\end{enumerate}

In addition to that we also removed the scaling parameter from being a parameter of the algorithm
\begin{enumerate}[label=\textcolor{seccolor}\textbullet]
	\item the scaling parameter $s$ is strictly defined at every step, not a degree of freedom of the algorithm.
\end{enumerate}

Summarizing this section, we can present the final set of instructions for RG smoothing:
\ttfamily
\begin{Figure}
\begin{parchment}[algorithm]
\begin{enumerate}[label=\textcolor{myred}\theenumi]
	\item enumerate data points;
	\item build a chain from tangent vectors (the direction is not important);
	\item find optimal scaling parameter $s$ using eq. \ref{eq:s};
	\item rescale the chain, i.e. do similar transformations as in eq.~1 and eq.~2 but for chosen~$s$;
	\item repeat steps (3) and (4) if needed.
\end{enumerate}
\end{parchment}
\end{Figure}
\rmfamily

The only freedom left (after fixing the chain) is the number of iteration steps.
So the number of repetitions (item 5 in the algorithm) is the only tuning parameter of the method.
\section*{Example}
In this part I will show how the algorithm works for an arbitrary scatter plot example.

\begin{figure*}[t]
	\vspace{-1.5cm}
	\hspace{-1.cm}
	\includegraphics[width=1.1\textwidth,trim={0.0cm 0.0cm 0.0cm 0.0cm}, clip]{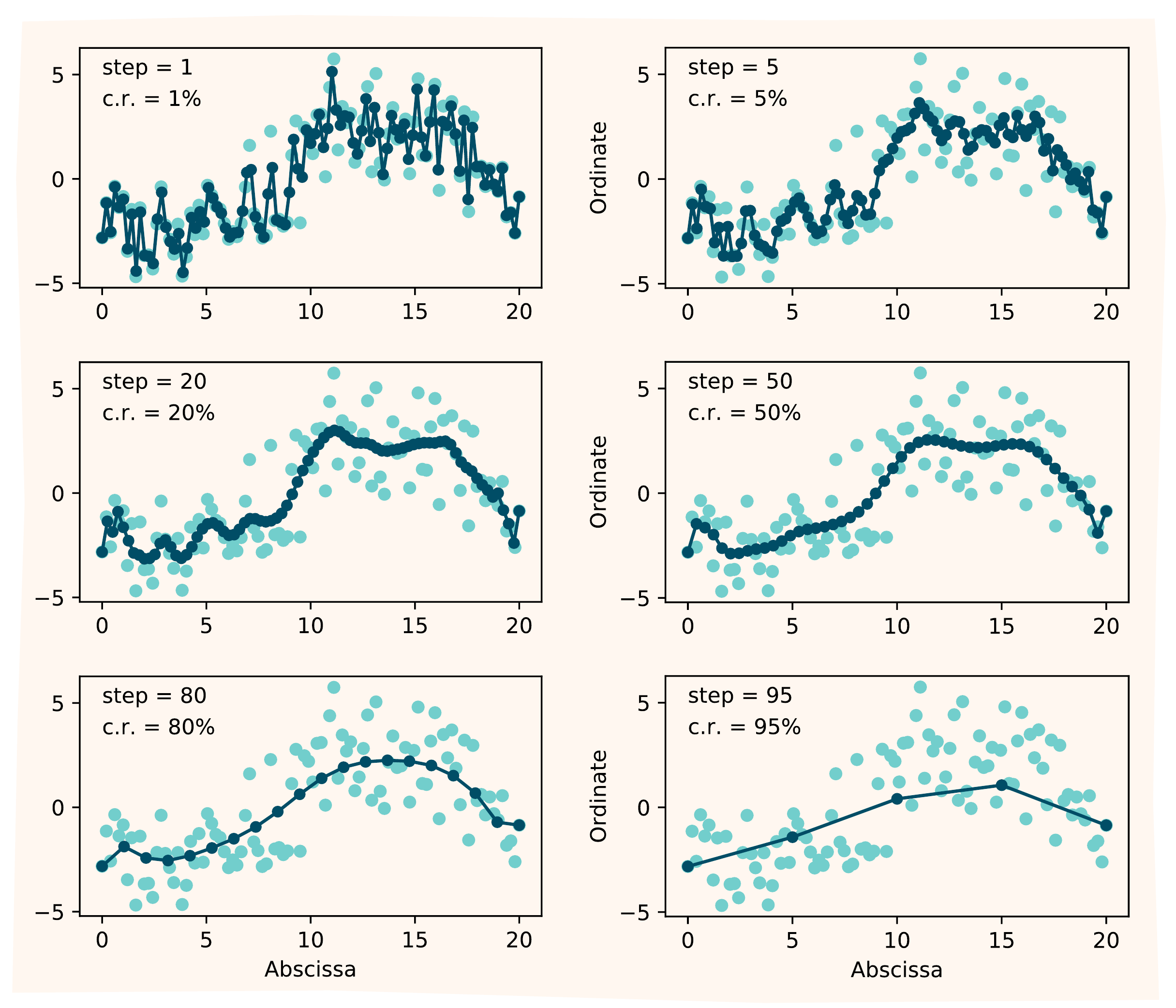}
	\vspace{-1.0cm}
	\caption{Smoothing example.}
\label{fig:example}
\end{figure*}

I generated 101 data points with regular x-grid.
The result of smoothing with RG method in shown on Figure~\ref{fig:example} with dark curves.
``Step" is a number of iterations we did.
So for the first plot we applied our algorithm only once.

The reduction of the number of data points can be expressed as a \hl{compression ratio} which is marked as 
``c.r." in the Figure and defined as:
\begin{equation*}
c.r. =\left( 1- \frac{N_p+1}{N+1} \right)  \times 100\%,
\end{equation*}
where $N+1$ is the number of data points in the original data, and $N_p+1$ is the number of points after the $p$-th step of smoothing.

If you take a look at the plots with 50 and 80 steps in the figure you can notice the different treatment of the end points.
Since we fix the first and the last point, they are included to the final chain without any corrections.
In other words the errors in the first and the last points remain untouched.

Another feature which should be noticed can be seen in the last two figures: for 80 and 95 steps.
In the last plot the $x$ coordinates for our new smoothed points are: 0, 5, 10, 15, 20.
It means that \hl{the property of equidistanty of the original data points is remained}.
Indeed, if we look at our transformations eq.~(1, 2) then it is obvious that we can write those transformation for $x$ components ans $y$ components independently.
The algorithm perform rescaling with parameter $s$ in all directions along basis vectors independently.

\section*{Conclusion}
A novel algorithm for smoothing was presented.
The idea behind the method is quite simple, however it is based on solid mathematical formulation of  renormalization group theory.
From that theory the algorithm got not only the name (RG smoothing), but very important features: data compression and recursive nature.

It leads to lots of advantages of RG smoothing algorithm and here I only list some of them together with possible implementations:

\begin{enumerate}[label=\textcolor{myred}\textbullet]
	\item can be used for flow-data analysis, where the question of data reduction is crucial;
	\item can be used for trajectory smoothing, because it can handle non-regular grids and self crossings;
	\item does not depend on the origin and characteristics of the noise;
	\item the algorithm is identical for any number of dimensions;
	\item can be used for image contour compression;
	\item can be implemented directly in hardware, since the algorithm is iterative and each step is simple;
	\item has only one tuning parameter - the number of iterative steps;
	\item the algorithm is non-parametric: number of parameters does not grow with thr size of the system;
	\item if the grid was regular initially then this property is conserved;
	\item can be used for contour recognition.
\end{enumerate}

The downside which still should be improved is how to handle the ends.





\section*{for Feedback}
This is a preprint that is work in progress.
Since my own background is in physics, I am interested in feedback and potential collaboration on the article itself with people that are knowledgable in fields such as signal processing, pattern recognition, etc.

If you have any ideas about possible applications and improvements you are very welcome to contact me.

\hlb{E-mail: anna.sinelnikova@physics.uu.se}

\end{multicols}

\vspace{1cm}
\noindent This work was supported by the Swedish research council and the Knut and Alice Wallenberg foundation through the Wallenberg Academy Fellow grant of J. Nilsson.
\end{document}